\documentclass[10pt,journal]{IEEEtran}
\usepackage{graphicx}          
\usepackage{amsmath}           
\usepackage{amssymb}
\usepackage{algorithm}
\usepackage{algpseudocode}
\usepackage{hyperref}
\usepackage[latin1]{inputenc}  

\usepackage{enumerate}
\usepackage{epstopdf}
\usepackage{color}

\usepackage{amsthm}
\newtheorem{thm}{Theorem}

\theoremstyle{definition}

\theoremstyle{remark}

\theoremstyle{plain}
\newtheorem{lem}[thm]{Lemma}

\DeclareMathOperator{\E}{E}
\DeclareMathOperator*{\argmin}{arg\,min}

\newcommand{\mbs}[1]{\boldsymbol{#1}}
\newcommand{\what}[1]{\widehat{#1}}

\newcommand{\wtilde}[1]{\widetilde{#1}}

\begin{document}

\title{Weighted SPICE: A Unifying Approach for Hyperparameter-Free Sparse Estimation}
\author{Petre Stoica, Dave Zachariah, Jian Li\thanks{ Corresponding
    author: dave.zachariah@it.uu.se. Ph: +46739513234. This work was supported in part by the National Science Foundation under Grant No. CCF-1218388, the Office of Naval Research (ONR) under Grant No. N00014-12-1-0381, and the U.S. Army Research Laboratory and the U.S. Army Research Office under contract/grant number No. W911NF-11-2-0039. The views and conclusions contained herein are those of the authors and should not be interpreted as necessarily representing the official policies or endorsements, either expressed or implied, of the U.S. Government. The U.S. Government is authorized to reproduce and distribute reprints for Governmental purposes notwithstanding any copyright thereon.}}


\maketitle

\begin{abstract}
In this paper we present the SPICE approach for sparse
parameter estimation in a framework that unifies it with other
hyperparameter-free methods, namely LIKES, SLIM and IAA.\footnote{All abbreviations used in this paper are explained at the end of the Introduction.} Specifically,
we show how
the latter methods can be interpreted as variants of an adaptively
reweighted SPICE method. Furthermore, we establish a connection
between SPICE and the $\ell_1$-penalized LAD estimator as well as the square-root
LASSO method. We evaluate the four
methods mentioned above in a generic sparse regression problem and in
an array processing application.
\end{abstract}

\section{Introduction}

During the past two decades, sparse parameter estimation for the ubiquitous
linear model
\begin{equation}
\mbs{y} = \mbs{B}\mbs{x} + \mbs{e} \in \mathbb{C}^N, \quad \mbs{B} \in \mathbb{C}^{N \times M}
\label{eq:linmodel}
\end{equation}
has become an important problem in signal processing, statistics and
machine learning
\cite{HastieEtAl2009_statisticallearning,Candes&Wakin2008_introcompressed,Elad2010_sparse}, with
applications ranging from spectral analysis and direction-of-arrival
estimation to magnetic resonance imaging and biomedical analysis
\cite{MalioutovEtAl2005, BourguignonEtAl2007, LustigEtAl2008,
  Tibshirani1996}. In this model $\mbs{x} \in \mathbb{C}^M$
denotes the unknown
sparse parameter vector and $\mbs{y} \in \mathbb{C}^N$ is the vector
of observations with typically $M \gg N$. The matrix of regressors $\mbs{B}$ is assumed to be
given and the unknown noise $\mbs{e}$ is assumed to have zero mean. 
\textcolor{black}{For $M > N$, the problem is ill-posed, unless the
  knowledge about the sparsity of $\mbs{x}$ is exploited.} 

\textcolor{black}{Exploiting sparsity also
  enables one to tackle nonlinear estimation problems. Consider, for example, a nonlinear model which consists of the superposition of
  an unknown number of mode vectors
$$\mbs{y} = \sum_{i} \mbs{h}(\tilde{\mbs{\theta}}_i) \tilde{x}_i + \mbs{e},$$ 
where $\mbs{h}(\tilde{\mbs{\theta}})$ is a given function of unknown parameters
$\tilde{\mbs{\theta}}_i \in \Theta$. Each mode has an amplitude
$\tilde{x}_i$. This model is relevant to spectral analysis and related
applications. By griding the parameter space $\Theta$ using a sufficiently large number of points $\{
\mbs{\theta}_k \}^M_{k=1}$, we can approximate the
nonlinear model by means of a dictionary of mode vectors evaluated at the fixed grid points $\mbs{B}=[\mbs{h}(\mbs{\theta}_1) \:
\cdots \: \mbs{h}(\mbs{\theta}_M)]$ such that $\mbs{x} = [x_1 \: \cdots \:
x_M]^\top$ becomes a sparse vector and \eqref{eq:linmodel}
applies. Identification of a nonzero component $x_k$
therefore concomitantly  identifies the corresponding gridpoint
$\mbs{\theta}_k$ which becomes an estimate of the nonlinear parameter.}
\textcolor{black}{Another example, that appears in applications of machine
  learning and statistics, is finding a general  input-output mapping
$$
y_i = h(\mbs{\theta}_i) + e_i,
$$ 
from data $\{ \mbs{\theta}_i, y_i \}^N_{i=1}$.  The  nonlinear mapping is modeled by a
sparse linear combination of kernel functions, $h(\mbs{\theta}_i) = \sum^N_{j=1} k(\mbs{\theta}_i,\mbs{\theta}_j)x_j$, where $\{ x_j\}$ denote the expansion
coefficients \cite{Roth2004_generalized}. Thus the sparse linear model \eqref{eq:linmodel} applies and identification of the
nonlinear mapping can be posed as a sparse parameter estimation
problem, where the $ij$th element of $\mbs{B}$
is given by the kernel function  $k(\mbs{\theta}_i, \mbs{\theta}_j)$.}

Many popular sparse estimation methods are
based on regularizing the least-squares method by penalizing a norm of the parameter
vector $\mbs{x}$, \textcolor{black}{in an attempt to strike a balance
  between data fidelity and parameter sparsity.} While such sparsifying methods can estimate $\mbs{x}$ in highly
underdetermined scenarios, most of them require the careful selection of user-defined regularization hyperparameters \cite{ChenEtAl2001_bpdn, Tropp2006_justrelax, Fuchs2004_sparserepresentation,EfronEtAl2004_lars}, cf. \cite{Malecki&Donoho2010_tunedcompressed} for a critical discussion.

Recently, a \textbf{sp}arse \textbf{i}terative \textbf{c}ovariance-based
\textbf{e}stimation method (\textsc{Spice}) was proposed which does
not require any hyperparameters, yet has good statistical properties \cite{StoicaEtAl2011_spicespectral,StoicaEtAl2011_spicearray,Stoica&Babu2012_spicelikes}. In this tutorial paper:
\begin{itemize}
\item We set out to derive four different hyperparameter-free methods in a
unifying \textsc{Spice}-like manner: the methods are \textsc{Spice},
\textsc{Likes},
\textsc{Slim} and \textsc{Iaa} \cite{Stoica&Babu2012_spicelikes,TanEtAl2011_slim,YardibiEtAl2010_iaa}. In the process we provide
insights into these methods, and derive new versions of each of them.
\item Furthermore, we establish the connection between
  \textsc{Spice} and $\ell_1$-penalized \textsc{Lad} as well as the square-root
  \textsc{Lasso} methods
  \cite{Wang2013_penLAD,BelloniEtAl2011_sqrtlasso}.

\item Finally, we evaluate the four methods in two
different scenarios: a generic sparse regression problem and a
direction-of-arrival estimation application.
\end{itemize}

\emph{Notation:} Matrices, vectors and scalars are distinguished by $\mbs{A}$, $\mbs{a}$ and $a$, respectively. $\text{Re}\{a\}$ denotes the real part of $a$. Defined variables are signified by $\triangleq$. $\mbs{A}^{1/2}$ is a matrix square-root of $\mbs{A}$ and $\mbs{A}^{-1/2}$ is its inverse. $\mbs{A} \otimes \mbs{B}$ and $\mbs{A} \odot \mbs{B}$ denote the Kronecker
and Khatri-Rao matrix products. $\mbs{A}^\top$ and
$\mbs{A}^*$ denote the transpose and Hermitian transpose of $\mbs{A}$. $\text{vec}(\mbs{A})$ is the column-wise vectorization of $\mbs{A}$. $\| \cdot \|_1$, $\| \cdot \|_2$ and $\| \cdot \|_F$ are
the $\ell_1$, $\ell_2$ and Frobenius norms, and $\|
\cdot \|_0$ denotes the $\ell_0$ `quasi-norm' which equals the
number of nonzero entries of a vector. $\text{tr}\{ \mbs{A}\}$ and
$|\mbs{A}|$ denote the trace and determinant of a square matrix
$\mbs{A}$. We use $\text{diag}(d_1, \dots, d_N)$ or
$\text{diag}(\mbs{d})$ to compactly denote a diagonal matrix with
entries from $\mbs{d}$. $\mbs{A} \succeq \mbs{B}$ signifies the Löwner order
between Hermitian matrices $\mbs{A}$ and $\mbs{B}$.  The Kronecker delta is
denoted by $\delta_{jk}$. The proper complex Gaussian distribution with mean $\mbs{\mu}$ and covariance matrix $\mbs{\Sigma}$ is denoted $\mathcal{CN}(\mbs{\mu}, \mbs{\Sigma})$. The probability of event $E$ is written as $\Pr\{ E \}$.

\emph{Abbreviations:} If and only if (iff). Subject to
(s.t.). With respect to (w.r.t.). Identically and
independently distributed (IID). Signal-to-noise ratio (SNR). Mean square error (MSE). Linear minimum mean square error
(LMMSE). Least squares (LS). Second-order cone program (SOCP).  Direction-of-arrival (DOA). Uniform linear array (ULA).  Least absolute deviation (\textsc{Lad}). Least absolute shrinkage and
selection operator (\textsc{Lasso}). Focal underdetermined system
solver (\textsc{Focuss}). Sparse iterative covariance-based estimation
(\textsc{Spice}). Likelihood-based estimation of sparse parameters
(\textsc{Likes}). Sparse learning via iterative minimization (\textsc{Slim}). Iterative adaptive approach (\textsc{Iaa}).

\section{Brief review of the basic Spice approach}

\textsc{Spice} was introduced as a covariance fitting approach in
\cite{StoicaEtAl2011_spicespectral,StoicaEtAl2011_spicearray,Stoica&Babu2012_spicelikes}. In
what follows we consider the one-snapshot case of \eqref{eq:linmodel},
but the method is easily extended to the multisnapshot case as we show
in Appendix~\ref{app:multisnapshot}.
Consider the following `model' for the covariance matrix of the data vector $\mbs{y}$:
\begin{equation}
\begin{split}
\mbs{R} &= \mbs{B}
\begin{bmatrix}
p_1 &        & 0\\
   & \ddots & \\
  0  &        & p_M
\end{bmatrix}
\mbs{B}^* + \begin{bmatrix}
p_{M+1} &        & 0 \\
   & \ddots & \\
  0  &        & p_{M+N}
\end{bmatrix} \\
&= \mbs{A}\mbs{P} \mbs{A}^* \in \mathbb{C}^{N \times N},
\end{split}
\label{eq:covariancemodel}
\end{equation}
where
\begin{equation*}
\mbs{A} \triangleq [\mbs{B} \; \mbs{I}_N], \quad
\mbs{P} \triangleq \text{diag}(\mbs{p}),
\end{equation*}
and where $\mbs{A} = [\mbs{a}_1 \: \mbs{a}_2 \: \cdots \: \mbs{a}_{M+N}] \in
\mathbb{C}^{N \times (M+N)}$ and $\mbs{p}  \triangleq [p_1, \dots,
p_{M+N}]^\top \in \mathbb{R}^{M+N}_+$. The covariance matrix
$\mbs{R}(\mbs{p})$ is a function of the parameters $\{ p_k \}$ which can
be interpreted as the variances of $\{ x_k \}$ and $\{ e_k \}$.
In the next section, we will discuss the covariance model
\eqref{eq:covariancemodel} in more detail. While it appears to assume
that $\{ x_k, e_k \}$ are uncorrelated, this should not be interpreted as a
restriction, as will be explained.

In the spectral analysis applications of sparse parameter estimation, the
main goal is to estimate $\{ p_k \}$. In most of the other
applications, the goal is the estimation of $\mbs{x}$. Even in the
latter case, there exists a class of methods (that includes those
discussed here) which first obtain estimates $\{ \hat{p}_k \}$ of $\{
p_k \}$ and then, if desired, estimate $\mbs{x}$ via the usual LMMSE estimator
formula \cite{Kay1998_estimationtheory}:
\begin{equation}
\hat{x}_k = \hat{p}_k \mbs{a}^*_k \what{\mbs{R}}^{-1} \mbs{y}, \: k=1,\dots,M,
\label{eq:MMSEMAP}
\end{equation}
where $\what{\mbs{R}} = \mbs{A}\what{\mbs{P}}\mbs{A}^*$. As we
show in the next sections, this estimate also occurs naturally within
an augmented version of the \textsc{Spice} approach.
An alternative is to use the Capon formula \cite{Stoica&Moses2005}:
\begin{equation}
\hat{x}_k = \frac{\mbs{a}^*_k \what{\mbs{R}}^{-1} \mbs{y} }{
  \mbs{a}^*_k \what{\mbs{R}}^{-1} \mbs{a}_k },  \: k=1,\dots,M.
\label{eq:WLSBLUE}
\end{equation}
In general one can expected \eqref{eq:WLSBLUE} to be a less biased estimate than \eqref{eq:MMSEMAP}, but \eqref{eq:MMSEMAP} to have a smaller MSE. Interestingly, if the same $\what{\mbs{P}}$ is used in both \eqref{eq:MMSEMAP} and \eqref{eq:WLSBLUE} then:
\begin{equation}
|\hat{x}_k|_{(3)} \leq |\hat{x}_k|_{(4)}.
\end{equation}
In particular this means that the estimate \eqref{eq:MMSEMAP} of $\mbs{x}$ is always sparser than \eqref{eq:WLSBLUE}. This fact follows from the following simple result:
\begin{lem}
\begin{equation}
\hat{p}_k \leq \frac{1}{\mbs{a}^*_k \what{\mbs{R}}^{-1} \mbs{a}_k}. \:
\label{eq:lemmaPinv}
\end{equation}
\end{lem}
The proof of this lemma, as well as of the subsequent ones, can be found in
Appendix~\ref{app:proofs}.

If a $K$-sparse estimate of $\mbs{x}$ is
desired, that is an estimate $\{ \hat{x}_k \}$ where only $K$ elements are
nonzero, then we can apply the LS method to \eqref{eq:linmodel} where we
retain only the columns of $\mbs{B}$ whose indices
correspond to the $K$ largest peaks of $\{ \hat{p}_k \}^M_{k=1}$.

\textsc{Spice} estimates $\{ p_k \}$ by minimizing the following covariance fitting criterion:
\begin{equation}
\begin{split}
\| \mbs{R}^{-1/2}(\mbs{y}\mbs{y}^* - \mbs{R}) \|^2_F &= \text{tr}\{ (\mbs{y}\mbs{y}^* - \mbs{R}) \mbs{R}^{-1}(\mbs{y}\mbs{y}^* - \mbs{R}) \} \\
&=\| \mbs{y} \|^2_2 \mbs{y}^* \mbs{R}^{-1} \mbs{y} + \text{tr}\{ \mbs{R}
\} +  \text{const.},
\label{eq:SPICEcovariancematch}
\end{split}
\end{equation}
or equivalently,
\begin{equation}
\mbs{y}^* \mbs{R}^{-1} \mbs{y} + \frac{1}{\| \mbs{y} \|^2_2} \sum^{M+N}_{k=1} w_k p_k, \quad w_k = \| \mbs{a}_k \|^2_2.
\label{eq:SPICEunnormalized}
\end{equation}
Next we note the following result:
\begin{lem}
Let
$$\hat{\mbs{p}} = \argmin_{\mbs{p}} \: g(\mbs{p}), \; g(\mbs{p}) = \mbs{y}^* \mbs{R}^{-1}(\mbs{p}) \mbs{y} + c^2 \sum^{M+N}_{k=1} w_k p_k,$$
where $\: c > 0$, and let
$$\hat{\bar{\mbs{p}}} = \argmin_{\bar{\mbs{p}}} \: f(\bar{\mbs{p}}), \; f(\bar{\mbs{p}}) = \mbs{y}^* \mbs{R}^{-1}(\bar{\mbs{p}}) \mbs{y} + \sum^{M+N}_{k=1} w_k \bar{p}_k.$$
Then
$$\hat{\bar{\mbs{p}}} = c \hat{\mbs{p}}.$$
\end{lem}
Note also that a uniform scaling of $\{ \hat{p}_k \}$ leaves $\{ \hat{x}_k \}$
unchanged whether using LMMSE, Capon or LS.
It follows from these observations that the constant factor $\|
\mbs{y} \|^2_2$ in \eqref{eq:SPICEunnormalized} can be omitted. Thus
we can reformulate the \textsc{Spice} criterion as:
\begin{equation}
\min_{\{p_k \}} \; \mbs{y}^* \mbs{R}^{-1}\mbs{y} + \sum^{M+N}_{k=1} w_k p_k.
\label{eq:SPICE}
\end{equation}
When $\| \mbs{a}_k \|_2 \equiv \text{const}.$, the weights in \eqref{eq:SPICE} can be replaced by 1's.

The problem in \eqref{eq:SPICE} is convex, namely an SOCP \cite{StoicaEtAl2011_spicearray}, and hence it can be solved globally \cite{Boyd&Vandenberghe2004_convex}. Rather than
solving it by an off-the-shelves SOCP code, the following iterative cyclic minimizer, which monotonically decreases \eqref{eq:SPICE} at each iteration and converges globally \cite{StoicaEtAl2011_spicespectral,StoicaEtAl2011_spicearray}, was found to be preferable from a computational standpoint:
\begin{equation}
\hat{p}^{i+1}_k = \hat{p}^{i}_k | \mbs{a}^*_k \what{\mbs{R}}^{-1}_i
\mbs{y} | / \sqrt{w_k}, \quad (\text{\textsc{Spice}}_\text{a})
\label{eq:SPICEiteration}
\end{equation}
where $\quad k=1,2,\dots, M+N$ and $\what{\mbs{R}}_i
= \mbs{A}\text{diag}(\hat{\mbs{p}}^i)\mbs{A}^*$ denotes the
covariance matrix estimate at iteration $i$; we use a subindex `a' for the \textsc{Spice} algorithm
in \eqref{eq:SPICEiteration} to differentiate it from a variation that
will be presented later on, see Section~\ref{sec:reweightedspice}
below.


We remark on the fact that we have allowed the noise $\mbs{e}$ to have
different powers, say $\sigma^2_1 = p_{M+1}$, $\sigma^2_2 =
p_{M+2}, \dots, \sigma^2_N = p_{M+N} $, in different data samples for three reasons:
\begin{itemize}
\item notational simplicity (it allows treating the noise powers
  similarly to $\{ p_k \}^M_{k=1}$ and not differently as is the case
  when the condition $\sigma^2_k \equiv \sigma^2, \forall k$, is enforced).

\item generality (in some applications, $\sigma^2_1,
  \sigma^2_2, \dots, \sigma^2_N$,  may well be different from one another).

\item even if the noise powers are uniform, $\sigma^2_k \equiv \sigma^2$,
  $\forall k$, but we either do not know that or prefer not to impose
  this condition, \textsc{Spice} with different noise powers works
  well. Indeed, the degradation of accuracy compared with that
  achieved by imposing $\sigma^2_k \equiv \sigma^2$, $\forall k$, is
  not significant -- we explain why this is so in the next section. On
  the other hand, if we are sure that $\sigma^2_k \equiv \sigma^2$,
  $\forall k$, and want to enforce this condition, then we can do so
  with only some minor modifications of the algorithms (see
  \cite{StoicaEtAl2011_spicespectral,StoicaEtAl2011_spicearray,Stoica&Babu2012_spicelikes}
  and Appendix~\ref{app:identicalnoise}).
\end{itemize}
Finally, we note that the form of \textsc{Spice}$_\text{a}$ iteration,
\eqref{eq:SPICEiteration}, is similar to that associated with
\textsc{Focuss} \cite{Gorodnitsky&Rao1997_focuss}; the main difference
between the two methods lies in the way the noise powers are
treated: \textsc{Focuss} assumes that the noise powers are identical
\emph{and} given (possibly estimated by some other method), whereas
\textsc{Spice} does not make this restrictive assumption. 

\section{On the covariance model and the link of Spice to
  $\ell_1$-penalized Lad}
\label{sec:covmodelandL1LAD}

There are several important questions about the covariance model in \eqref{eq:covariancemodel}:
\begin{enumerate}[a)]
\item Assume that $\mbs{x}$ and $\mbs{e}$ are drawn from correlated
  distributions (i.e. distributions whose covariance matrices are not
  diagonal). Then will \textsc{Spice}, and the other estimation
  methods discussed later, still work despite seemingly relying on the diagonal covariance matrix in \eqref{eq:covariancemodel}? Note that in the Bayesian approach to sparse estimation (see e.g. \cite{Tipping2001_sblrvm}), \eqref{eq:covariancemodel} is viewed as a `prior information' -- however this does not offer any satisfactory answer to the above equation, as priors are not `forgotten' in problems with many more unknowns than data samples (i.e. $M \gg N)$, as in the case considered here.

\item In fact what do $\{ \hat{p}_k \}$ estimate? Do they estimate $\{|x_k|^2 \}$?

\item If indeed $\sigma^2_k \equiv \sigma^2$, $\forall k$, do we
  significantly degrade the accuracy by using a nonuniform noise
  power model as in \eqref{eq:covariancemodel}?

\item Is \eqref{eq:covariancemodel} a unique description, i.e. for a given $\mbs{P}$ can we find a $\bar{\mbs{P}} \neq \mbs{P}$ such that $\mbs{A}\mbs{P}\mbs{A}^* = \mbs{A}\bar{\mbs{P}}\mbs{A}^*$?
\end{enumerate}
We will provide answers to questions a)-c) by establishing the
connection between the \textsc{Spice} criterion in \eqref{eq:SPICE}
and the $\ell_1$-penalized \textsc{Lad} criterion. Then we will address the
question in d) by means of a separate analysis.

To understand the role of $\mbs{p}$ in the estimation of $\mbs{x}$, we rewrite the criterion in \eqref{eq:SPICE} in terms of the original model \eqref{eq:linmodel}, with the help of the following result.
\begin{lem}
Let
$$\mbs{S} = \begin{bmatrix} p_{M+1} & & 0 \\ & \ddots &  \\ 0 & & p_{M+N}\end{bmatrix}, \; \mbs{\Pi} = \begin{bmatrix} p_{1} & & 0 \\ & \ddots &  \\ 0 & & p_{M}\end{bmatrix}.$$
Then
\begin{equation}
\min_{\mbs{x}} \: (\mbs{y} - \mbs{B}\mbs{x})^* \mbs{S}^{-1} (\mbs{y} - \mbs{B}\mbs{x}) + \sum^{M}_{k=1} |x_k|^2 /p_k = \mbs{y}^* \mbs{R}^{-1} \mbs{y}
\label{eq:SPICEalt}
\end{equation}
and the minimum value occurs at
\begin{equation}
\hat{\mbs{x}} = \mbs{\Pi}\mbs{B}^* \mbs{R}^{-1} \mbs{y}.
\label{eq:SPICEaltsolution1}
\end{equation}
\end{lem}

It follows from the above lemma that the minimizer of the \textsc{Spice} criterion in \eqref{eq:SPICE} can also be obtained by minimizing the following function (w.r.t. \emph{both} $\mbs{x}$ \emph{and} $\mbs{p}$):
\begin{equation}
\sum^N_{k=1} | y_k - \mbs{b}^*_k \mbs{x} |^2/p_{M+k} + \sum^M_{k=1} |x_k|^2/p_k \; + \sum^{M+N}_{k=1} w_k p_k,
\label{eq:SPICEalt2}
\end{equation}
where $\mbs{b}^*_k$ denotes the $k$th row of $\mbs{B}$.
Minimization of \eqref{eq:SPICEalt2} w.r.t. $\{ p_k \}$ yields:
\begin{equation}
\begin{split}
p_k &= |x_k|/\sqrt{w_k}, \; k=1,\dots, M \\
p_{M+k} &= |y_k - \mbs{b}^*_k \mbs{x}|/\sqrt{w_{M+k}}, \; k=1,\dots, N.
\end{split}
\label{eq:SPICEalt2_pk1}
\end{equation}
Insertion of \eqref{eq:SPICEalt2_pk1} into \eqref{eq:SPICEalt2} gives:
\begin{equation}
\| \mbs{W}^{1/2}_1( \mbs{y} - \mbs{B} \mbs{x}) \|_1 + \left \|
\mbs{W}^{1/2}_2 \mbs{x}
\right \|_1,
\label{eq:L1LAD}
\end{equation}
where $\mbs{W}_1 = \text{diag}(w_{M+1}, \dots, w_{M+N})$ and
$\mbs{W}_2 = \text{diag}(w_{1}, \dots, w_{M})$; this is recognized as a (weighted) $\ell_1$-penalized \textsc{Lad} criterion \cite{Wang2013_penLAD}.

The above analysis has several implications, some for \textsc{Lad}:
\begin{itemize}
\item The $\ell_1$-penalized \textsc{Lad} estimate of $\mbs{x}$ can be obtained using the \textsc{Spice}$_{\text{a}}$ algorithm, \eqref{eq:SPICEiteration}, to estimate $\{ p_k \}$ and then get $\hat{\mbs{x}}$ from
    \eqref{eq:SPICEaltsolution1} (note that \eqref{eq:SPICEaltsolution1} is identical to \eqref{eq:MMSEMAP}). For the complex-valued data case, \textsc{Spice} can be expected to be faster than other convex programming techniques that are used to get $\hat{\mbs{x}}$ from \eqref{eq:L1LAD} directly.

\item If the condition $\sigma^2_1 = \sigma^2_2 = \cdots = \sigma^2_N$
  is enforced, then the \textsc{Spice} approach was shown in
  \cite{Babu&StoicaEtAl2013_spicesqrtlasso,RojasEtAl2013_notespice} to be equivalent to the
  square-root \textsc{Lasso} method of
  \cite{BelloniEtAl2011_sqrtlasso} (see also Appendix~\ref{app:identicalnoise} for a more
  direct proof of this equivalance result). This fact establishes an interesting connection between square-root \textsc{Lasso} and $\ell_1$-penalized \textsc{Lad}.
\end{itemize}
and some for \textsc{Spice}:
\begin{itemize}
\item The \textsc{Spice} estimates $\{ \hat{p}_k \}$ of $\{ p_k \}$
  are not estimates of $|x_k|^2$ and $|e_k|^2$ but of (scaled versions
  of) the square-roots of these quantities. However, when $w_k$ is an estimate
  of $1/p_k$, or a quantity related to $1/p_k$ (see Section~\ref{sec:reweightedspice}), then it follows from
  \eqref{eq:SPICEalt2_pk1} that $\{ \hat{p}_k \}$ estimates $|x_k|^2$.


\item \textsc{Spice} will still work
  even if the $\{ x_k \}$ and $\{ e_k \}$ in \eqref{eq:linmodel} are
  drawn from correlated distributions; indeed, when \textsc{Spice} is
  viewed from the perspective of its equivalence with
  $\ell_1$-penalized \textsc{Lad} (or square-root \textsc{Lasso}), its performance
  does not depend significantly on the way $\mbs{x}$ and $\mbs{e}$
  were generated because the performance of the $\ell_1$-penalized \textsc{Lad}
  or square-root \textsc{Lasso} is not strongly dependent on that
  \cite{Wang2013_penLAD,BelloniEtAl2011_sqrtlasso}; in this
  light, $\{ p_k \}$ and the `covariance model' in
  \eqref{eq:covariancemodel} can be viewed as being nothing but
  instruments employed to achieve the equivalence proven
  above; and hence not for necessarily providing a true description of
  the data covariance matrix. Similarly, by not imposing the condition
  $\sigma^2_1 = \sigma^2_2 = \cdots = \sigma^2_N$, when this was known
  to be true, we basically do nothing but choose to use
  $\ell_1$-penalized \textsc{Lad} in lieu of square-root \textsc{Lasso}, and
  the difference in accuracy between the latter methods is usually
  not significant.

\end{itemize}

In the above discussion we have provided answers to questions
a)-c). Next, we turn our attention to question d). The
parameterization/description \eqref{eq:covariancemodel} of $\mbs{R}$
is unique iff there is no diagonal matrix $\mbs{Q} =
\text{diag}(\mbs{q}) \neq \mbs{0}$ (where $\mbs{q} \in \mathbb{R}^{M+N}$)  which is such that:
\begin{equation}
\mbs{A}\mbs{Q}\mbs{A}^* = \mbs{0}
\label{eq:AQA0}
\end{equation}
and
\begin{equation}
\mbs{P} + \mbs{Q} \succeq \mbs{0},
\label{eq:PpQ0}
\end{equation}
Equation \eqref{eq:AQA0} can be re-written as:
\begin{equation*}
\begin{split}
\text{vec}(\mbs{A}\mbs{Q}\mbs{A}^*) &= \text{vec}\left(\sum^{M+N}_{k=1} q_k \mbs{a}_k \mbs{a}^*_k \right) \\
&= [\mbs{a}^{*\top}_1 \otimes \mbs{a}_1 \: \cdots \:  \mbs{a}^{*\top}_{M+N} \otimes \mbs{a}_{M+N}] \mbs{q} = \mbs{0},
\end{split}
\end{equation*}
i.e., equivalently,
\begin{equation}
(\mbs{A}^{*\top} \odot \mbs{A}) \mbs{q} = \mbs{0}.
\label{eq:KhatriRaonullspace}
\end{equation}
Also, for easy reference, we can write \eqref{eq:PpQ0} as $\mbs{p} +\mbs{q} \geq \mbs{0}$.
In the analysis of \eqref{eq:KhatriRaonullspace}, the rank of the matrix $\wtilde{\mbs{A}} \triangleq \mbs{A}^{*\top} \odot \mbs{A}$ clearly is an essential factor:
\begin{enumerate}[i)]
\item Assume that $M < N$ (i.e. $M+N < 2N$) and that any $N$ columns of $\mbs{A}$ are linearly independent. Then $\text{rank}(\wtilde{\mbs{A}}) = M+N$ \cite{SidiropoulosEtAl2000_parfac} and thus \eqref{eq:KhatriRaonullspace} implies $\mbs{q} = \mbs{0}$, which means that the description \eqref{eq:covariancemodel} of $\mbs{R}$ is unique.

\item If $M+N < N^2$, then $\wtilde{\mbs{A}}$ has full rank (equal to $M+N$) for almost any matrix $\mbs{A}$ (assumed to be drawn from a continuous distribution) \cite{JiangEtAl2001_multidimharmonic} and hence \eqref{eq:covariancemodel} is generically unique.

\item For $M+N > N^2$, $\text{rank}(\wtilde{\mbs{A}}) < M+N$ and there
  exists $\mbs{q} \neq \mbs{0}$ that satisfy
  \eqref{eq:KhatriRaonullspace}. In this scenario one must consider
  two cases. Let $r$ be the maximum integer such that any $r$ columns of
  $\tilde{\mbs{A}}$ are linearly independent. Then, if $\| \mbs{p} \|_0
  > r$,
  \eqref{eq:covariancemodel} is not unique, whereas if $\| \mbs{p} \|_0
  \leq r$ then \eqref{eq:covariancemodel} may be unique or nonunique
  depending on the instance of $\mbs{A}$ and $\mbs{p}$ under consideration.
To see this, let $\bar{\mbs{A}}$
  denote the matrix made from the columns of $\tilde{\mbs{A}}$ corresponding to
  the nonzero elements of $\mbs{p}$. When  $\| \mbs{p} \|_0
  > r$ there exists
  a vector $\bar{\mbs{q}} \neq \mbs{0}$ such that $\bar{\mbs{A}}
  \bar{\mbs{q}} = \mbs{0}$. By appending $\bar{\mbs{q}}$ with zeros we
  can therefore form a vector $\mbs{q} $ that fulfills both
  $\tilde{\mbs{A}}\varepsilon \mbs{q} = \mbs{0}$ and $\mbs{p} +
  \varepsilon \mbs{q} \geq 0$, for a sufficiently small $\varepsilon$. On the other hand, when  $\| \mbs{p} \|_0
  \leq r$ then such a vector $\mbs{q}$ in the nullspace of
  $\tilde{\mbs{A}}$ may or may not satisfy \eqref{eq:PpQ0} depending
  on whether the signs of the coefficients which do not belong to the
  support of $\mbs{p}$ are all the same.


\end{enumerate}

\section{Likes, Slim, Iaa (and new versions thereof) as (re)weighted
  Spice}
\label{sec:reweightedspice}

Consider the \textsc{Spice} fitting criterion in \eqref{eq:SPICE} with
general weights $\{ w_k > 0 \}$ (possibly different than the
weights in \eqref{eq:SPICE}). For fixed weights, \eqref{eq:SPICE} is a
convex function of $\{ p_k \}$, which can be globally minimized, for
example, by the algorithm in \eqref{eq:SPICEiteration}. In the
following we will derive \eqref{eq:SPICEiteration} by using a
\emph{gradient approach} that is simpler than the cyclic minimization
approach employed in
\cite{StoicaEtAl2011_spicespectral,StoicaEtAl2011_spicearray,Stoica&Babu2012_spicelikes}. The
gradient approach is also more flexible in that it suggests
alternatives to \eqref{eq:SPICEiteration} which may be interesting in their own right.

\subsection{Weighted \textsc{Spice}}

The derivative of \eqref{eq:SPICE} w.r.t. $p_k$ is equal to
\begin{equation}
- \mbs{y}^* \mbs{R}^{-1} \frac{\partial \mbs{R}}{\partial p_k} \mbs{R}^{-1} \mbs{y} + w_k = -|\mbs{a}^*_k \mbs{R}^{-1} \mbs{y}|^2 + w_k.
\end{equation}
Consequently, the $(i+1)$th iteration of a gradient algorithm (with variable step length) applied to \eqref{eq:SPICE} is given by:
\begin{equation}
\hat{p}^{i+1}_k = \hat{p}^{i}_k - \rho^i_k( w_k - |\mbs{a}^*_k \what{\mbs{R}}^{-1}_i \mbs{y}|^2 ),
\label{eq:p_kgradientstep}
\end{equation}
where $\what{\mbs{R}}_i$ is made from $\{ \hat{p}^i_k \}$, as before, and the step
size $\rho^i_k$ must be non-negative
\begin{equation}
\rho^i_k \geq 0.
\label{eq:rhopositive}
\end{equation}
Because $\{ p_k \geq 0 \}$ by definition, we shall also choose $\rho^i_k$ such that:
\begin{equation}
\hat{p}^i_k \geq 0 \: \Rightarrow \: \hat{p}^{i+1}_k \geq 0.
\label{eq:p_kpositive}
\end{equation}
Let us choose
\begin{equation}
\rho^i_k = \frac{\hat{p}^i_k}{w_k + w^{1/2}_k|\mbs{a}^*_k \what{\mbs{R}}^{-1}_i \mbs{y}|}
\label{eq:rhochoice1}
\end{equation}
which satisfies \eqref{eq:rhopositive}. A simple calculation gives:
\begin{equation*}
\begin{split}
\hat{p}^{i+1}_k &= \frac{\hat{p}^i_k w_k + \hat{p}^i_k w^{1/2}_k |\mbs{a}^*_k \what{\mbs{R}}^{-1}_i \mbs{y}| - \hat{p}^i_k w_k + \hat{p}^i_k |\mbs{a}^*_k \what{\mbs{R}}^{-1}_i \mbs{y}|^2 }{w_k + w^{1/2}_k |\mbs{a}^*_k  \what{\mbs{R}}^{-1}_i \mbs{y}|} \\
&= \frac{\hat{p}^i_k|\mbs{a}^*_k \what{\mbs{R}}^{-1}_i \mbs{y}|}{w^{1/2}_k},
\end{split}
\end{equation*}
that is,
\begin{equation}
\hat{p}^{i+1}_k = \hat{p}^i_k |\mbs{a}^*_k \what{\mbs{R}}^{-1}_i
\mbs{y}| / w^{1/2}_k \quad (\text{\textsc{Spice}}_\text{a})
\label{eq:SPICEiterationalt1}
\end{equation}
and thus \eqref{eq:p_kpositive} is satisfied too. Note that when $w_k = \| \mbs{a}_k \|^2_2$,
\eqref{eq:SPICEiterationalt1} is nothing but the \textsc{Spice}$_\text{a}$
algorithm in \textcolor{black}{equation} \eqref{eq:SPICEiteration}, whose derivation above is
more direct than the derivation in
\cite{StoicaEtAl2011_spicespectral,StoicaEtAl2011_spicearray,Stoica&Babu2012_spicelikes}
which was based on cyclically minimizing an augmented criterion
function. 

As already mentioned, the gradient approach is also more flexible in the sense that $\rho^i_k$ in \eqref{eq:p_kgradientstep} can be chosen in several different ways than \eqref{eq:rhochoice1} to obtain alternative algorithms to \eqref{eq:SPICEiterationalt1}. A particularly simple such choice (that satisfies \eqref{eq:rhopositive}) is:
\begin{equation}
\rho^i_k = \hat{p}^i_k/w_k
\label{eq:rhochoice2}
\end{equation}
which leads to
\begin{equation}
\hat{p}^{i+1}_k = \hat{p}^i_k |\mbs{a}^*_k \what{\mbs{R}}^{-1}_i
\mbs{y}|^2 / w_k   \quad (\text{\textsc{Spice}}_\text{b})
\label{eq:SPICEiterationalt2}
\end{equation}
(therefore \eqref{eq:p_kpositive} is satisfied as well). When $w_k = \| \mbs{a}_k \|^2_2$,
\eqref{eq:SPICEiterationalt2} minimizes the same criterion as
\eqref{eq:SPICEiterationalt1} and will therefore be referred to as
\textsc{Spice}$_{\text{b}}$. Both algorithms share the same stationary
points, but they may have
different rates of convergence. In particular observe that the step length in \eqref{eq:rhochoice1}
is smaller than \eqref{eq:rhochoice2},  when both are evaluated using the same $\{ \hat{p}^i_k \}$.

In the next sections we will consider different choices of the weights
than \textsc{Spice}'s, which will lead to other hyperparameter-free
methods, namely \textsc{Likes},
\textsc{Slim} and \textsc{Iaa}. Unlike \textsc{Spice}, whose weights are
constant, these algorithms use data-dependent weights that change with
the iteration. 

\subsection{\textsc{Likes}}

The current problem of estimating $\{ p_k \}$ from $\mbs{y}$ is not a
standard one especially owing to the fact that $M+N$ = number of
unknowns $\gg 2N$ = number of (real-valued) data. Even so, the
analysis in \cite{OtterstenEtAl_comet}, as well as data-whitening
considerations, suggest that a possibly (statistically) better
covariance matching criterion than \eqref{eq:SPICEcovariancematch} is
the following one:
\begin{equation}
\| \mbs{R}^{-1/2} (\mbs{y} \mbs{y}^* - \mbs{R})\what{\mbs{R}}^{-1/2} \|^2_F,
\label{eq:covfittingalt}
\end{equation}
where $\what{\mbs{R}}$ is an available estimate of $\mbs{R}$. A
straightforward calculation shows that \eqref{eq:covfittingalt} can be re-written as:
\begin{equation}
\begin{split}
&\text{tr}\left[ (\mbs{y}\mbs{y}^* - \mbs{R}) \what{\mbs{R}}^{-1} (\mbs{y}\mbs{y}^* - \mbs{R}) \mbs{R}^{-1} \right] \\&= (\mbs{y}^* \what{\mbs{R}}^{-1} \mbs{y})(\mbs{y}^* {\mbs{R}}^{-1} \mbs{y}) + \text{tr}(\what{\mbs{R}}^{-1} \mbs{R}) + \text{const.} \\
&=(\mbs{y}^* \what{\mbs{R}}^{-1} \mbs{y})(\mbs{y}^* {\mbs{R}}^{-1} \mbs{y}) + \sum^{M+N}_{k=1} (\mbs{a}^*_k \what{\mbs{R}}^{-1} \mbs{a}_k)p_k + \text{const.}
\end{split}
\label{eq:tracecovfittingalt}
\end{equation}

In view of Lemma~2 we can omit the constant factor $(\mbs{y}^* \what{\mbs{R}}^{-1} \mbs{y})$ in \eqref{eq:tracecovfittingalt}, which leads to the following weighted \textsc{Spice} criterion:
\begin{equation}
\mbs{y}^* {\mbs{R}}^{-1} \mbs{y} + \sum^{M+N}_{k=1} w_k p_k, \quad w_k = (\mbs{a}^*_k \what{\mbs{R}}^{-1} \mbs{a}_k).
\label{eq:SPICEalt3}
\end{equation}
Unlike \textsc{Spice}'s weights, which are data independent, the $\{
w_k \}$ in \eqref{eq:SPICEalt3} depend on the data (via
$\what{\mbs{R}}$). Note that $w_k$ in \eqref{eq:SPICEalt3} can be
interpreted as the Capon estimate of $1/p_k$ (see
e.g. \cite{Stoica&Moses2005}). This means that the penalty term in
\eqref{eq:SPICEalt3} is an approximation of $\| \mbs{p} \|_0$ rather
than just being proportional to the $\ell_1$-norm of $\mbs{p}$ as for
\textsc{Spice}. It is well known that the $\ell_0$-(quasi)norm is the most sensible measure of the sparsity of a parameter vector because it is not dependent on the size of the elements of that vector, as is the $\ell_1$-norm (see, e.g., \cite{CandesEtAl2008_reweightedmin} for a general discussion on this aspect).

It follows from the above discussion that the weights in \eqref{eq:SPICEalt3} are intuitively a more appealing choice than the \textsc{Spice}'s weights in \eqref{eq:SPICE}. The data-dependent weights in \eqref{eq:SPICEalt3} can be updated in the following way:
\begin{enumerate}[i)]
\item Fix $\what{\mbs{R}}$ in $\{ w_k \}$ and use \eqref{eq:SPICEiterationalt1} or \eqref{eq:SPICEiterationalt2} to minimize \eqref{eq:SPICEalt3}, or at least monotonically decrease this function for a pre-specified number of iterations.

\item Update $\what{\mbs{R}}$ in \eqref{eq:SPICEalt3}, and the weights $\{ w_k \}$, and go to step i).
\end{enumerate}
This leads to the following iterative schemes:
\begin{equation}
\hat{p}^{i+1}_k = \hat{p}^{i}_k |\mbs{a}^*_k \what{\mbs{R}}^{-1}_i \mbs{y}| / ( \mbs{a}^*_k \what{\mbs{R}}^{-1}_\ell \mbs{a}_k )^{1/2}, \quad (\text{\textsc{Likes}}_{\text{a}})
\label{eq:LIKES}
\end{equation}
or, alternatively,
\begin{equation}
\hat{p}^{i+1}_k = \hat{p}^{i}_k |\mbs{a}^*_k \what{\mbs{R}}^{-1}_i \mbs{y}|^2 / ( \mbs{a}^*_k \what{\mbs{R}}^{-1}_\ell \mbs{a}_k ) \quad (\text{\textsc{Likes}}_{\text{b}}).
\label{eq:LIKESprime}
\end{equation}
Initially, we set $\ell =
0$, and the above updates are executed as follows:
\begin{enumerate}
\item Iterate for $i = \ell, \ell+1, \dots, \ell+m-1$, where $m$ is the number
  of iterations in which the weights are kept fixed.

\item Reset $\ell \leftarrow \ell + m$, and go to 1).
\end{enumerate}

The algorithm in \eqref{eq:LIKES} is recognized as \textsc{Likes} \cite{Stoica&Babu2012_spicelikes},
whereas the one in \eqref{eq:LIKESprime} is a new version. To
distinguish between them we have designated them as \textsc{Likes}$_{\text{a}}$ and
\textsc{Likes}$_{\text{b}}$, respectively. Because these algorithms
update the weights in \eqref{eq:SPICEalt3}, they can only be
interpreted as minimizers of the criterion in \eqref{eq:SPICEalt3}
using the weights obtained \emph{at convergence}. This does not say
much as to the convergence properties of \eqref{eq:LIKES} or
\eqref{eq:LIKESprime}, an aspect that will be addressed in the next
section. Here we only note that the two iterative algorithms above clearly have the same stationary points. However, their rates of convergence to a stationary point may be different from one another.

\subsection{\textsc{Slim}}

Consider \eqref{eq:SPICEalt3} with different weights:
\begin{equation}
w_k = 1/\hat{p}_k.
\label{eq:SLIMweights}
\end{equation}
The corresponding penalty term in \eqref{eq:SPICEalt3} would then be a
more direct approximation of $\| \mbs{p} \|_0$ than when $w_k =
\mbs{a}^*_k \what{\mbs{R}}^{-1}\mbs{a}_k$ as for \textsc{Likes}. In
fact it follows from Lemma~1 that the weights in
\eqref{eq:SLIMweights} are \emph{larger} than \textsc{Likes}' weights. Consequently the use of
\eqref{eq:SLIMweights} should yield sparser estimates of $\{ p_k \}$
than \textsc{Likes} does. Note that this interpretation is valid as long as the weights are kept fixed and therefore it does not extend necessarily to the case in which the weights are updated (because in the latter case different weights lead to different estimates of $\{ p_k \}$ and hence the weights at different iterations do not correspond to the same $\{ \hat{p}_k \}$ any longer). However, empirical evidence suggests that the above observation remains typically valid even in that case.

Using \eqref{eq:SLIMweights} in \eqref{eq:SPICEiterationalt1} and \eqref{eq:SPICEiterationalt2} yields the algorithms:
\begin{equation}
\hat{p}^{i+1}_k = (\hat{p}^{i}_k)^{3/2} |\mbs{a}^*_k
\what{\mbs{R}}^{-1}_i \mbs{y}| \quad (\text{\textsc{Slim}}_\text{a})
\label{eq:SLIM}
\end{equation}
and
\begin{equation}
\hat{p}^{i+1}_k = (\hat{p}^{i}_k)^2 |\mbs{a}^*_k \what{\mbs{R}}^{-1}_i
\mbs{y}|^2 \quad (\text{\textsc{Slim}}_\text{b})
\label{eq:SLIMprime}
\end{equation}
where \eqref{eq:SLIMprime} is recognized as \textsc{Slim}
\cite{TanEtAl2011_slim} (more precisely, an extension of the
\textsc{Slim}-0 algorithm in the cited reference to the case of
different noise powers) and \eqref{eq:SLIM} is a new version
thereof that we call \textsc{Slim}$_\text{a}$. Most comments made in the previous subsection about the \textsc{Likes} algorithm apply to \eqref{eq:SLIM} and \eqref{eq:SLIMprime} as well. In particular, \eqref{eq:SLIM} and \eqref{eq:SLIMprime} clearly have the same stationary points.

\subsection{\textsc{Iaa}}

The weights in \eqref{eq:SLIMweights} were larger than \textsc{Likes}'. Next consider the following weights:
\begin{equation}
w_k = \hat{p}_k (\mbs{a}^*_k \what{\mbs{R}}^{-1} \mbs{a}_k)^2
\label{eq:IAAweights}
\end{equation}
which, in view of Lemma~1, are \emph{smaller} than \textsc{Likes}
weights (whenever both sets of weights are computed from the same $\{
\hat{p}_k \}$). The estimates of $\{ p_k \}$ corresponding to
\eqref{eq:IAAweights} can therefore be expect to be \emph{less sparse}
than \textsc{Likes} estimates; and this fact, despite the cautionary
note following \eqref{eq:SLIMweights}, is confirmed by empirical evidence.

Using \eqref{eq:IAAweights} in \eqref{eq:SPICEiterationalt1} and \eqref{eq:SPICEiterationalt2}, in the same fashion as done above for \textsc{Slim}, we get:
\begin{equation}
\hat{p}^{i+1}_k = (\hat{p}^i_k)^{1/2} |\mbs{a}^*_k \what{\mbs{R}}^{-1}_i \mbs{y}| / (\mbs{a}^*_k \what{\mbs{R}}^{-1}_i \mbs{a}_k) \quad (\text{\textsc{Iaa}}_\text{a})
\label{eq:IAA}
\end{equation}
and
\begin{equation}
\hat{p}^{i+1}_k = |\mbs{a}^*_k \what{\mbs{R}}^{-1}_i \mbs{y}|^2 / (\mbs{a}^*_k \what{\mbs{R}}^{-1}_i \mbs{a}_k)^2  \quad (\text{\textsc{Iaa}}_\text{b})
\label{eq:IAAprime}
\end{equation}

The same comments, made previously on the \textsc{Likes} and \textsc{Slim}
algorithms, apply verbatim to \textsc{Iaa}$_{\text{a}}$ and
\textsc{Iaa}$_{\text{b}}$ as well. Note that \textsc{Iaa}$_\text{b}$
concides with the original \textsc{Iaa} algorithm introduced in
\cite{YardibiEtAl2010_iaa} whereas \textsc{Iaa}$_\text{a}$ is a new version.


\section{Statistical interpretations and convergence properties}
\label{sec:statisticalinterpretations}

\subsection{\textsc{Spice}}

The \textsc{Spice} algorithms minimize the convex covariance fitting
criterion in \eqref{eq:SPICE}, and they can be shown to be
globally convergent from any initial estimate $\{ \hat{p}_k > 0 \}$ (\cite{StoicaEtAl2011_spicespectral,StoicaEtAl2011_spicearray,Stoica&Babu2012_spicelikes}). This
property basically follows from the convexity of the problem, and the
fact that both \textsc{Spice}$_\text{a}$ and \textsc{Spice}$_\text{b}$
monotonically decrease the optimization criterion (as explained in
Appendix \ref{app:cyclic}).

The other algorithms discussed here also globally minimize their
corresponding covariance fitting criteria provided that the weights
are kept fixed. This is a useful property as long as the weights are
reasonable approximations of $1/p_k$. However, when the weights are
continuously updated, as in \eqref{eq:SPICEalt3}, \eqref{eq:SLIMweights} and \eqref{eq:IAAweights}, this property is no longer valid and a separate analysis is needed to provide statistical interpretations of these algorithms, as well as analyze their convergence properties, see the next subsections.

\subsection{\textsc{Likes}}

Under the covariance model in \eqref{eq:covariancemodel} and the additional Gaussian data assumption, the negative log-likelihood function of $\mbs{y}$ is (to within an additive constant):
\begin{equation}
\mbs{y}^* \mbs{R}^{-1} \mbs{y} + \ln|\mbs{R}| .
\label{eq:negativeloglikelihood}
\end{equation}
The first term in \eqref{eq:negativeloglikelihood} is a convex
function whereas the second is an increasing concave function of $\{ p_k > 0 \}$ \cite{Stoica&Babu2012_spicelikes}. This implies that the second term in \eqref{eq:negativeloglikelihood} acts as a sparsity-inducing penalty. The previous fact also means that the function $\ln|\mbs{R}|$ in \eqref{eq:negativeloglikelihood} is majorized by its tangent plane at any point $\hat{\mbs{p}}$, that is by the following linear function of $\{ p_k \}$ (after omitting some uninteresting additive constants):
\begin{equation}
\begin{split}
\sum^{M+N}_{k=1} \frac{\partial ( \ln |\mbs{R}| )}{ \partial p_k }
\Bigl|_{p_k = \hat{p}_k} p_k &=
\sum^{M+N}_{k=1} \text{tr}\left[ \mbs{R}^{-1} \frac{\partial \mbs{R}
  }{ \partial p_k } \right] \Bigl|_{p_k = \hat{p}_k} p_k \\
 &= \sum^{M+N}_{k=1} (\mbs{a}^*_k \what{\mbs{R}}^{-1} \mbs{a}_k)p_k
\label{eq:majorization}
\end{split}
\end{equation}
Inserting \eqref{eq:majorization} into
\eqref{eq:negativeloglikelihood} we get the criterion in
\eqref{eq:SPICEalt3}. The \textsc{Likes} algorithms decrease
\eqref{eq:SPICEalt3} at each iteration  (see, once again Appendix~\ref{app:cyclic}) and therefore, by the
properties of majorization-minimization approaches (e.g. \cite{Stoica&Selen2004_cyclicminimization}), they decrease
\eqref{eq:negativeloglikelihood} monotonically. This fact implies that
the sequence of \textsc{Likes} estimates converges to a local minimum
of \eqref{eq:negativeloglikelihood}, or at least that it contains such
a convergent sub-sequence \cite{Zangwill1969_nonlinear}. Because the
current estimation problem is not a standard one, as already
mentioned, convergence to a minimum of the negative
log-likelihood function in \eqref{eq:negativeloglikelihood} does not
automatically guarantee good statistical properties; nevertheless it
is an interesting statistical interpretation of \textsc{Likes}.

\emph{Remark}: In the light of the above discussion, the
\textsc{Spice} criterion can also be related to
\eqref{eq:negativeloglikelihood} by replacing the penalty term
$\ln|\mbs{R}|$ in \eqref{eq:negativeloglikelihood} with $\text{tr}\{ \mbs{R} \}$. The criterion associated with \textsc{Slim} can be similarly interpreted, see below.

\subsection{\textsc{Slim}}

If $\ln | \mbs{R}|$ in \eqref{eq:negativeloglikelihood} is replaced by
$\ln|\mbs{P}|$, which is also an increasing concave function of $\{ p_k
> 0 \}$ and thus can serve as a penalty term, we obtain the criterion:
\begin{equation}
\mbs{y}^* \mbs{R}^{-1} \mbs{y} + \sum^{M+N}_{k=1} \ln p_k.
\label{eq:negativelogliklihoodalt}
\end{equation}
The tangent plane for the second term in \eqref{eq:negativelogliklihoodalt}, at any  $\{ \hat{p}_k\}$, is given by (to within an additive constant);
\begin{equation}
\sum^{M+N}_{k=1} \frac{\partial (\ln p_k)}{\partial p_k} \Bigr|_{p_k =
  \hat{p}_k} p_k = \sum^{M+N}_{k=1} \frac{1}{\hat{p}_k} p_k.
\label{eq:majorizationalt}
\end{equation}
Insertion of \eqref{eq:majorizationalt} in
\eqref{eq:negativelogliklihoodalt} yields a majorizing function for
\eqref{eq:negativelogliklihoodalt} that coincides with the
\textsc{Slim} criterion \eqref{eq:SLIMweights}. Consequently,
similarly to what was concluded following \eqref{eq:majorization}
about \textsc{Likes}, the \textsc{Slim} algorithms generate a sequence
of estimates that monotonically decreases
\eqref{eq:negativelogliklihoodalt} and converges to a minimum of
this function, or at least comprises a sub-sequence that does
so.


\subsection{\textsc{Iaa}}

Both \textsc{Likes} and \textsc{Slim}
monotonically decrease a cost function of the
form
\begin{equation}
\mbs{y}^* \mbs{R}^{-1}\mbs{y} + h(\mbs{p}),
\label{eq:costfunctionalform}
\end{equation}
where $h(\mbs{p})$ is an increasing concave function. For \textsc{Iaa}, on the
other hand, no function of this form can be found. 

The proof is by contradiction. If
indeed a concave
function $h(\mbs{p})$ existed whose derivatives w.r.t. $p_k$ were the
weights $\{ w_k \}$ of \textsc{Iaa}, i.e.,
$$ \frac{\partial h(\mbs{p})}{\partial p_k} =
p_k (\mbs{a}^*_k \mbs{R}^{-1} \mbs{a}_k)^2,$$ then the Hessian matrix
of that function would have the elements:
$$\frac{\partial^2 h(\mbs{p})}{ \partial p_k \partial p_j } =
 (\mbs{a}^*_k \mbs{R}^{-1}\mbs{a}_k)^2 \delta_{kj} - 2 p_k
 (\mbs{a}^*_k \mbs{R}^{-1}\mbs{a}_k)  |\mbs{a}^*_k \mbs{R}^{-1}\mbs{a}_j|^2.$$
But this matrix is not symmetric as required, let alone negative definite, and thus we reached a contradiction.

A partial statistical motivation of \textsc{Iaa} along with a local convergence
proof can be found in
\cite{YardibiEtAl2010_iaa,RobertsEtAl2010_iaalocalconv}. A more
definitive statistical interpretation of \textsc{Iaa} and a global analysis of
its convergence properties are open problems that await resolution. A
possible way of attacking these problems is to view the \textsc{Iaa} algorithms
as fixed-point iterations and attempt to make use of the available results
on the convergence of such iterations in the literature (see, e.g., \cite{Kelley1995_iterative}) to settle at least the question about \textsc{Iaa}'s convergence properties.

\subsection{Implementational aspects}
\emph{Version-a vs. version-b algorithms:} \label{sec:versionAvsB}
Empirical experience with the previous algorithms suggests that the convergence of
\textsc{Spice}$_\text{b}$ and \textsc{Likes}$_\text{b}$ can be significantly
slower than that of \textsc{Spice}$_\text{a}$ and
\textsc{Likes}$_\text{a}$. A plausible explanation for this follows
from the analysis in Appendix~\ref{app:cyclic}: when using \eqref{eq:SPICEiterationalt2} instead of
\eqref{eq:SPICEiterationalt1} we get equality in
\eqref{eq:gradcyc_inequalities}, instead of inequality, possibly
leading to a smaller reduction of the cost function.
Furthermore, the new \textsc{Iaa}$_\text{a}$, was found to work at least as well or better than the original \textsc{Iaa}$_\text{b}$. These findings suggest using
the a-versions of the algorithms rather than the b-versions.

\emph{Initialization and termination:} Unless otherwise stated, the algorithms are initialized with power estimates obtained from a matched filter,
$\hat{p}^0_k = | \mbs{a}^*_k \mbs{y} |^2 / \| \mbs{a}_k \|^4_2$, $\forall k$ and the
convergence tolerance for termination $\varepsilon$ in $\| \hat{\mbs{p}}^{i+1} - \hat{\mbs{p}}^{i} \|_2 / \|
\hat{\mbs{p}}^i \|_2 < \varepsilon$ is set to $10^{-3}$.  The
algorithms are set to terminate if the number of iterations exceeded
1000. 

\emph{\textsc{Spice}:} The implementation as in \eqref{eq:SPICEiterationalt1} follows the original setup of the algorithm \cite{StoicaEtAl2011_spicearray}, and was found to be numerically stable for all the tested cases.

\emph{\textsc{Likes}:} The original version of \textsc{Likes} was formulated as an iterative application of the \textsc{Spice} algorithm in which the weights are refined repeatedly \cite{Stoica&Babu2012_spicelikes}. \textsc{Likes} minimizes a nonconvex function with a number of local minima that typically increases as $N$ grows. Empirically we found that initializing the algorithm with the power estimates from \textsc{Spice}, as in the original formulation, produces better results than when using the matched filter. This is how we will initialize \textsc{Likes} in the numerical evaluations. Further, we update the weights as in \eqref{eq:LIKES} with $m=30$. It was found that too frequent updates led to performance degradation.

\emph{\textsc{Slim}:} As we have seen \textsc{Slim} decreases a cost function with a concave penality term. This function, however, lacks a global minimum; it assumes $-\infty$ if any power estimate is $0$. Therefore it is advisable to terminate after a small number of iterations, which is corroborated by empirical experience, cf. \cite{TanEtAl2011_slim}. Unlike \textsc{Spice}, which solves the powers of an $\ell_1$-penalized problem, \textsc{Slim} can be understood as a heuristic approach to approximate an $\ell_0$-penalized problem. We set the number of iterations, somewhat aribitrarily, to 5 in the numerical evaluations.

\emph{\textsc{Iaa}:} Empirically we found that when $N$ grows large 
numerical instabilities could occur due to numerical errors when computing $\mbs{a}^*_k \mbs{R}^{-1}_i \mbs{a}_k > 0$, which make the quantity complex-valued. We ensure that this quantity is real-valued when numerically evaluating the weights of \textsc{Iaa} and \textsc{Likes}, i.e., use $\text{Re}\{ \mbs{a}^*_k \mbs{R}^{-1}_i \mbs{a}_k \}$.

\emph{Remark:} In the interest of reproducible research we have made the codes for \textsc{Spice}, \textsc{Likes}, \textsc{Slim} and \textsc{Iaa}, as well as for the simulations in the subsequent section, available at \textcolor{black}{\texttt{https://www.it.uu.se/katalog/davza513}}.

\section{Numerical comparisons}

In this section we compare the four hyperparameter-free methods,
\textsc{Spice}, \textsc{Likes}, \textsc{Slim} and \textsc{Iaa}, by
means of numerical examples. \textcolor{black}{The standard \textsc{Lasso} with cross-validation based hyperparameter selection has already been compared with $\ell_1$-penalized \textsc{Lad} in
  \cite{Wang2013_penLAD}. In the cited paper and in \cite{BelloniEtAl2011_sqrtlasso}, the robustness of $\ell_1$-penalized \textsc{Lad} and square-root \textsc{Lasso} with respect to the hyperparameter choice was demonstrated and shown to be an important advantage over the standard \textsc{Lasso}}. Here, two different sparse parameter inference
problems are addressed for the linear model in \eqref{eq:linmodel} with
$\mbs{e} \sim \mathcal{CN}(\mbs{0}, \sigma^2 \mbs{I}_N)$. Note that despite generating noise with uniform powers,
we will not impose this constraint but rather use the general
algorithms derived in the previous sections.

First we
consider a generic regression problem with IID regressors, $b_{ij} \sim
\mathcal{CN}(0,1)$. In this case the cross-correlation
between the columns of $\mbs{B}$ is low. Next, we consider a
DOA estimation problem in which the adjacent columns of $\mbs{B}$
are highly correlated with each other. In both problems we let
$M=1000$.

We define the signal-to-noise ratio as
$\text{SNR} \triangleq \E[ \| \mbs{B} \mbs{x} \|^2_2 ] / \E[\| \mbs{e}
\|^2_2] = \sum_{k\in S} |x_k|^2 / \sigma^2$, where $S$ denotes the true
support set of nonzero coefficients. The performance metrics are evaluated using 1000 Monte Carlo
simulations. \textcolor{black}{We used a PC with Intel i7 3.4~GHz CPU and
  16~GB RAM. The algorithms were implemented in \textsc{Matlab} (MS Win7) in a rather direct manner without paying significant attention to computational details.}

\subsection{IID regressors}

The regressor matrix $\mbs{B}$ is randomized in each Monte Carlo
run. We consider $K$-sparse vectors $\mbs{x} \in \mathbb{C}^{1000}$,
where $K=3$, with a fixed support set $S = \{ 400, 420, 600
\}$. The nonzero coeffients $x_k = |x_k| e^{j \phi_k}$
have fixed powers, $\{ 1, 9, 4 \}$, respectively, and uniformly drawn
phases, for $k \in S$. The estimates $\hat{x}_k$ are computed using the
LMMSE formula \eqref{eq:MMSEMAP}. The Capon formula \eqref{eq:WLSBLUE}
produces less sparse estimates with higher MSE.

\begin{figure*}
  \begin{center}
    \includegraphics[width=2.0\columnwidth]{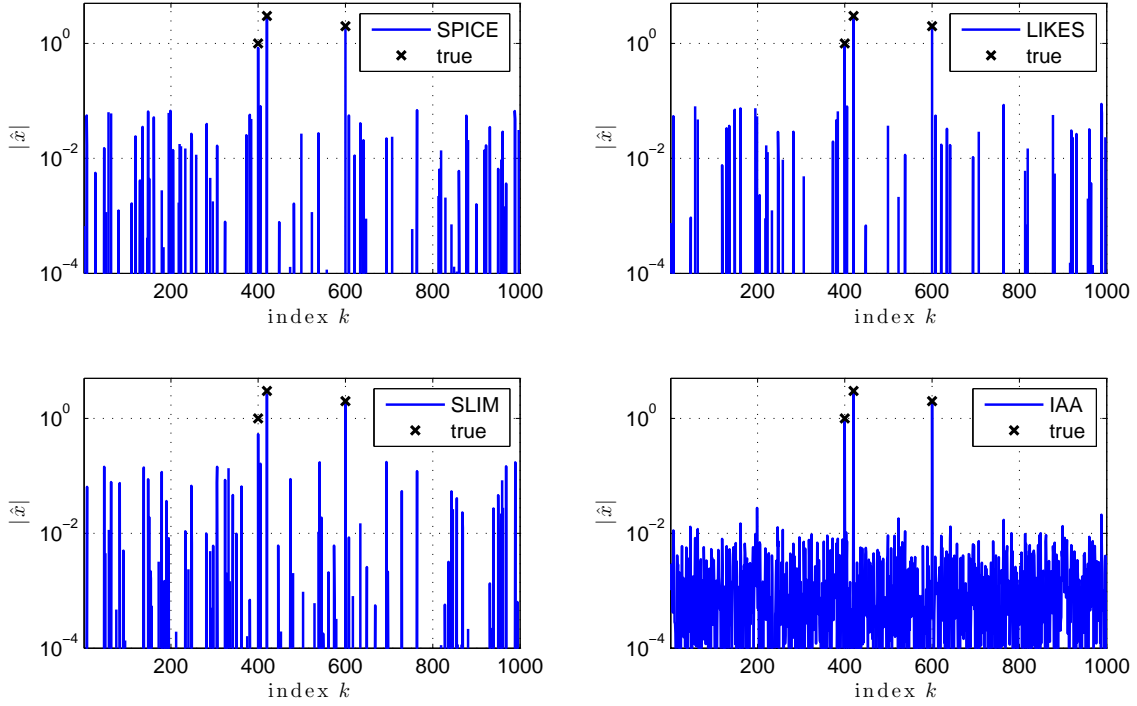}
  \end{center}
  \caption{Estimates $|\hat{x}_k|$ versus $k$ for a randomly selected
    realization. $N=35$ samples and $\text{SNR}=20$~dB.}
  \label{fig:IID_1snap}
\end{figure*}

Figure~\ref{fig:IID_1snap} illustrates the ability of the four
algorithms to locate the
active coefficients $\{ x_k \}_{k \in S}$ and provide reasonably small
estimates of $\{ x_k \}_{k \not \in S}$, for a randomly selected
realization. \textsc{Likes} and \textsc{Iaa} produce sparser respectively
denser estimates than \textsc{Spice}. Note that the magnitude of \textsc{Iaa} estimates
for $k \not \in S$ is substantially lower than for the other algorithms.
A plausible explanation of this is that the power estimates for $k
\not \in S$
capture a fraction of the residual power. Thus a `quasi-sparse' method
like \textsc{Iaa} will spread this residual power more evenly across $k \not
\in S$, than
a sparse method such as \textsc{Slim} which will concentrate it into
fewer nonzero estimates.

\begin{figure}
  \begin{center}
    \includegraphics[width=1.0\columnwidth]{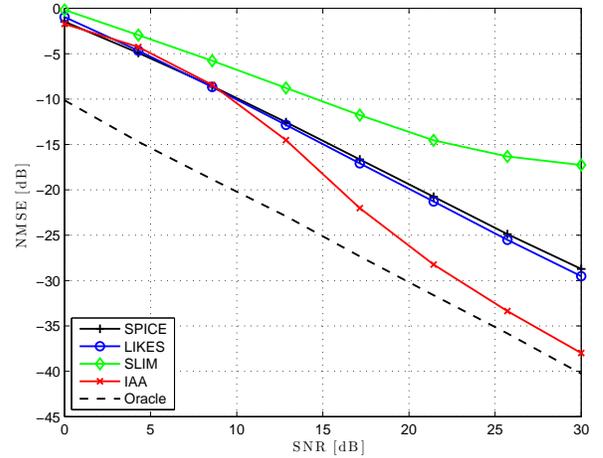}
  \end{center}
  \caption{Normalized MSE versus SNR for the IID regression problem, $N=35$ samples.}
  \label{fig:IID_MSEvsSNR}
\end{figure}

\begin{figure}
  \begin{center}
    \includegraphics[width=1.0\columnwidth]{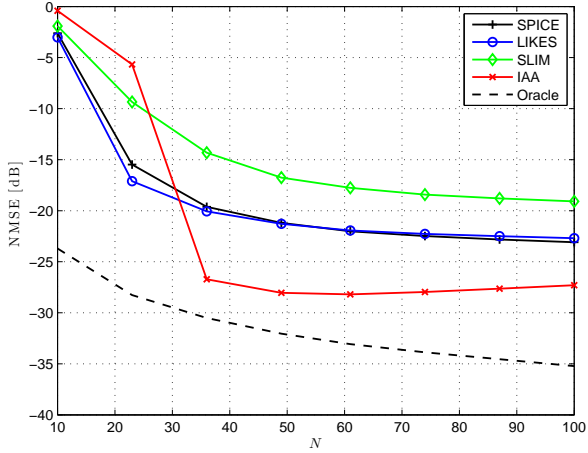}
  \end{center}
  \caption{Normalized MSE versus $N$ for the IID regression problem, $\text{SNR}=20$~dB.}
  \label{fig:IID_MSEvsN}
\end{figure}

\begin{figure}
  \begin{center}
    \includegraphics[width=1.0\columnwidth]{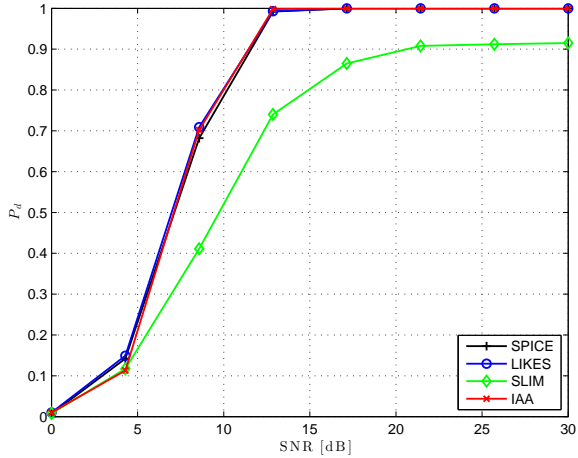}
  \end{center}
  \caption{Probability of correct support-set detection $P_d$ versus
    $\text{SNR}$ for the IID regression problem, $N=35$ samples.}
  \label{fig:IID_Pd_SNR}
\end{figure}

\begin{figure}
  \begin{center}
    \includegraphics[width=1.0\columnwidth]{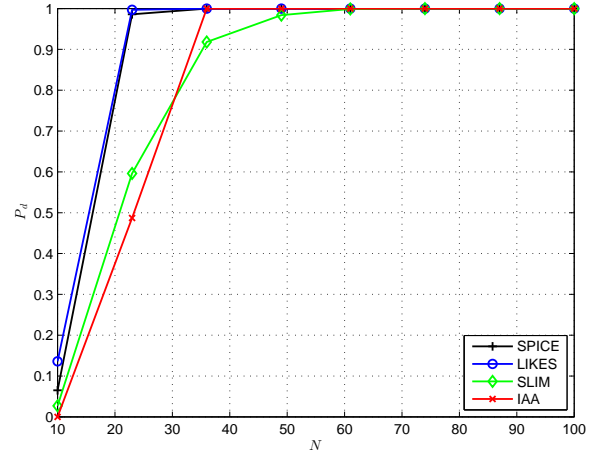}
  \end{center}
  \caption{Probability of correct support-set detection $P_d$ versus
    $N$ for the IID regression problem, $\text{SNR}=20$~dB.}
  \label{fig:IID_Pd_N}
\end{figure}

Figures~\ref{fig:IID_MSEvsSNR} and \ref{fig:IID_MSEvsN} show the mean
square error metric $\text{MSE}  \triangleq \E[ \|
\mbs{x} - \hat{\mbs{x}} \|^2_2 ]$, normalized by the signal power
$\E[\|\mbs{x}\|^2_2]$. This metric quantifies the ability of the
methods to localize $k \in S$ as well as provide reasonably small
estimates for $k \not \in S$. For reference we have added the performance of an
`oracle' estimator for which the unknown support set $S$ is
given; it computes the LS estimate for these coefficients, the
performance of which provides a lower MSE bound. Note that as
$M=1000$, the uniqueness condition $M+N < N^2$, is satisfied when $N
\geq 33$, cf. Section~\ref{sec:covmodelandL1LAD}. Observe that when
$N$ is above this threshold, \textsc{Iaa} performs better than the other
algorithms in terms of MSE. This
MSE reduction is mainly attributable to \textsc{Iaa}'s ability to provide
smaller coefficient estimates for $k \not \in S$.
The next two figures show plots of the support-set detection rate, $P_d \triangleq \Pr\{
S = \hat{S} \}$. We obtain the estimated support set, $\hat{S}$, for
each algorithm as the set of indices corresponding to the $K=3$ largest
values of $\hat{p}_k$, $k=1,\dots,M$. Figures~\ref{fig:IID_Pd_SNR}
and \ref{fig:IID_Pd_N} show $P_d$ as a function of $\text{SNR}$ and
$N$. We can see that $P_d$ approaches 1 for all algorithms as $N$
increases, and also that
\textsc{Spice} and \textsc{Likes} perform the best in the low sample scenario. The performance of the standard beamformer was too low for visibility and therefore omitted.

\begin{figure}
  \begin{center}
    \includegraphics[width=1.0\columnwidth]{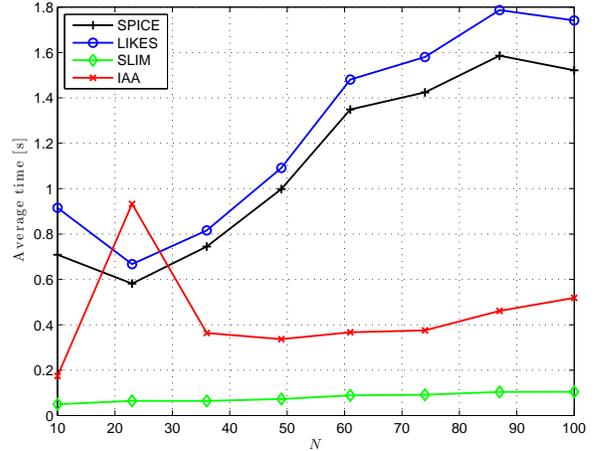}
  \end{center}
  \caption{Average computation time versus $N$ for the IID regression problem, $\text{SNR}=20$~dB.}
  \label{fig:IID_cputime}
\end{figure}

Finally, Figure~\ref{fig:IID_cputime} shows the average
computation time until convergence for each algorithm. While the
implementations are not carefully optimized, the figure should illustrate at least the
relative order of the algorithms. Noticeably, in the IID case with the
present signal dimensions, \textsc{Spice} tends to
be slower than \textsc{Slim} and \textsc{Iaa} which update
their weights adaptively. Not only does the performance of \textsc{Iaa} degrade when $N < 33$, but the algorithm tends to require more iterations until convergence. 

\subsection{Steering-vector regressors}

We now consider estimating the directions of arrival of the source
signals impinging on a uniform linear array
(ULA) with half-wavelength separation between elements. In this
problem the locations of the nonzero components of $\mbs{x}$ are of
interest rather than $\mbs{x}$ itself.
 The columns of $\mbs{B} = [\mbs{b}(\theta_1) \cdots \mbs{b}(\theta_M)]$
are given by the array
steering vector $\mbs{b}(\theta) = [1 \: e^{-j \kappa \sin\theta } \:
\cdots \: e^{-j (N-1) \kappa \sin\theta }]^\top$ \cite{Stoica&Moses2005}, and a uniform grid
of angles $\{ \theta_k \}^M_{k=1} \subset [-90^\circ, 90]$.\footnote{Here $\kappa = \omega_c d /c$, where $\omega_c$ is the signal frequency, $d$ is the element spacing and $c$ is the propagation velocity. We set $\kappa = \pi$.} We
consider $K=3$ sources
located at $\theta_k$,  $k \in S = \{ 400, 420, 600 \}$ on the
grid. This corresponds to DOAs at approximately $-18.1^\circ$,
$-14.5^\circ$ and $17.9^\circ$, respectively. As before the amplitudes
for $k \in S$ are generated as $x_k = |x_k| e^{j \phi_k}$
with fixed powers $\{ 1, 9, 4 \}$, respectively, and uniformly drawn
phases.

\begin{figure*}
  \begin{center}
    \includegraphics[width=2.0\columnwidth]{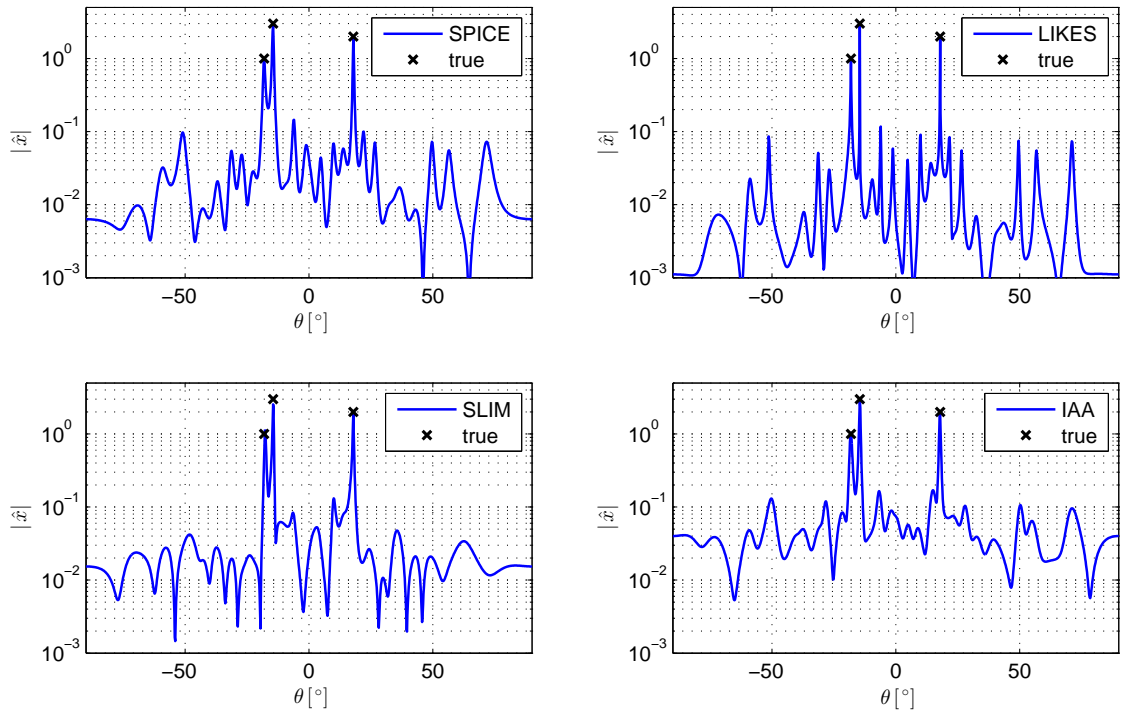}
  \end{center}
  \caption{Estimates $|\hat{x}_k|$ in  a randomly selected realization. $N=35$
    sensors and $\text{SNR}=20$~dB.}
  \label{fig:ULA_1snap}
\end{figure*}

\begin{figure}
  \begin{center}
    \includegraphics[width=1.0\columnwidth]{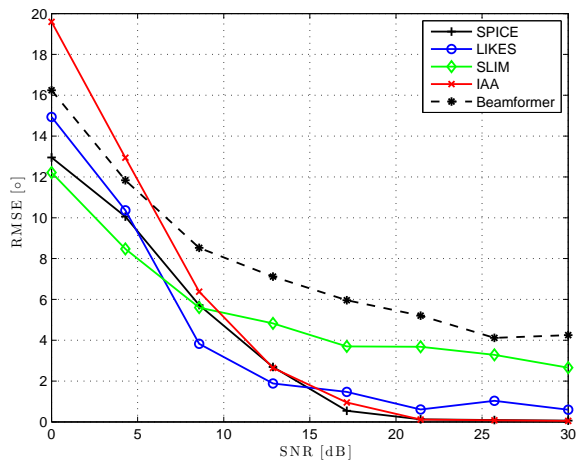}
  \end{center}
  \caption{Root mean square error of DOA estimates, per
    source, versus SNR. $N=35$ sensors.}
  \label{fig:ULA_RMSEvsSNR}
\end{figure}

\begin{figure}
  \begin{center}
    \includegraphics[width=1.0\columnwidth]{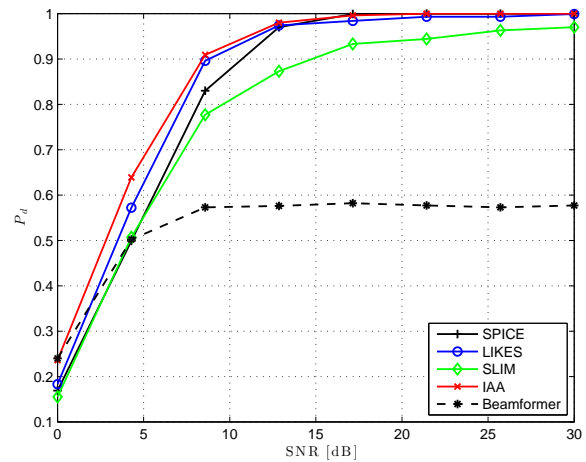}
  \end{center}
  \caption{Probability of detection $P_d$ versus SNR for
    steering-vector regressors, $N=35$ sensors.}
  \label{fig:ULA_PdvsSNR}
\end{figure}

Figure~\ref{fig:ULA_1snap} illustrates the ability of the four
algorithms to locate the sources and estimate their amplitudes in a
randomly selected
realization. The estimates $\hat{x}_k$ are computed using the
Capon formula \eqref{eq:WLSBLUE} which in the present case is less biased towards zero
than \eqref{eq:MMSEMAP}. Note that \textsc{Likes} produces sharper
spectral estimates than the other algorithms.
Next, we quantify the accuracy of the DOA estimates $\{ \hat{\theta}_k
\}$ obtained from the locations of the three peaks of $\{
\hat{p}_k \}^M_{k=1}$. In Figure \ref{fig:ULA_RMSEvsSNR} we plot the
root MSE per source, $\text{RMSE} \triangleq \sqrt{ \frac{1}{K}
  \E[ \| \mbs{\theta} - \hat{\mbs{\theta}} \|^2_2 ] }$, where
$\mbs{\theta}$ and $\hat{\mbs{\theta}}$ denote the vectors of ordered
DOAs and estimates, respectively. For reference, we have also included the
standard-beamformer performance. As SNR increases above 10~dB the
errors of \textsc{Spice}, \textsc{Likes} and \textsc{Iaa} fall well below the
RMSE of the beamformer. Figure \ref{fig:ULA_PdvsSNR}
shows the probability of detecting the $K$ sources within $\Delta
\theta$ degrees from the true DOA, $P_d \triangleq \Pr \{ |
\theta_i - \hat{\theta}_i | < \Delta \theta, \forall i \}$. Here we
set $\Delta \theta$ to half of the distance between the two closely-spaced
DOAs, i.e., $\Delta \theta = 1.8^\circ$. For this metric, \textsc{Iaa} turns
out to perform at least slightly better than the other
algorithms which all locate the peaks substantially better than the
beamformer. \textcolor{black}{For a further analysis of the resolution
  limit of sparse methods, see \cite{CandesEtAl2006_robust,FortunatiEtAl2014_csdoanalysis}.}

\begin{figure}
  \begin{center}
    \includegraphics[width=1.0\columnwidth]{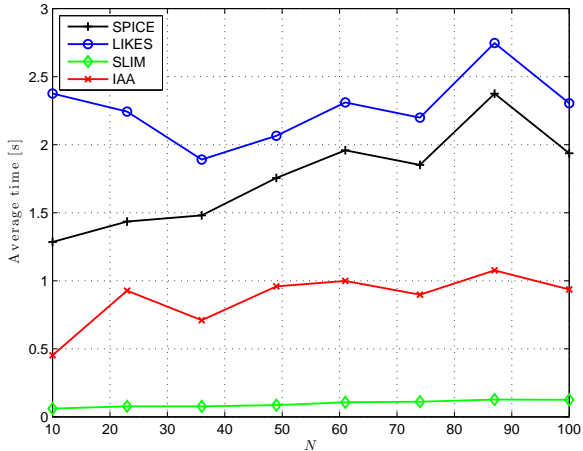}
  \end{center}
  \caption{Average computation time versus $N$ for steering-vector
    regressors, $\text{SNR}=20$~dB.}
  \label{fig:ULA_cputime}
\end{figure}

Figure~\ref{fig:ULA_cputime} illustrates the average computation time
versus $N$, and the order of the algorithms is the same as in Figure~\ref{fig:IID_cputime}. Recall that \textsc{Slim} is set to terminate after 5 iterations.








\section{Conclusions}

In this article we have presented a framework for sparse
parameter estimation based on the \textsc{Spice} fitting criterion and a gradient optimization approach. This framework was shown to have several appealing features:
\begin{itemize}

\item It unifies four hyperparameter-free methods, namely \textsc{Spice},
  \textsc{Likes}, \textsc{Slim} and \textsc{Iaa}, the latter three
  methods being instances of \textsc{Spice} with adaptive weights.

\item It enables further insights into the above four algorithms, including the analysis of their convergence properties and
  statistical performance. In particular, we showed
  how the weights used by these methods determine the sparsity of their corresponding estimates.

\item Finally, it makes it possible to derive new versions of the algorithms by
  considering different step-lengths in the gradient approach.
\end{itemize}
We also investigated the covariance model upon which the
\textsc{Spice} criterion is based, and:
\begin{itemize}

\item Provided identifiability conditions for this model.

\item  Showed that depending on whether the noise powers are modeled nonuniformly or
  uniformly, the \textsc{Spice} method coincides with the
  $\ell_1$-penalized \textsc{Lad} or the square-root \textsc{Lasso} problems. This fact also established a connection between the latter two methods.
\end{itemize}
The four hyperparameter-free methods were evaluated in two different
inference problems with IID and steering-vector regressors, respectively. The results indicated that:
\begin{itemize}
\item The a-version algorithms appear to be better than
  the b-versions in terms of convergence and statistical performance.

\item In problems with IID regressors both
  \textsc{Spice} and \textsc{Likes} perform similarly and they exhibit
  a graceful degradation as the number of samples decreases. For a
  sufficient number of samples, such that $M+N < N^2$, the
  `quasi-sparse' \textsc{Iaa} method, however, was found to \textcolor{black}{provide smaller
   parameter estimates} for the true zero coefficients.

\item In the steering-vector regressor case the peaks of the amplitude
  spectrum using the Capon formula were less biased towards zero than
  when using the LMMSE formula. \textsc{Likes} was
  computationally more demanding than the rest, but produced a sparser
  amplitude spectrum. For locating spectral peaks, however,
  \textsc{Iaa} was found to perform slightly better than the rest.
\end{itemize}

\appendices

\section{The multisnapshot case}
\label{app:multisnapshot}

The \textsc{Spice} criterion \eqref{eq:SPICEcovariancematch} extends to the multisnapshot
scenario as follows:
\begin{equation}
\begin{split}
\| \mbs{R}^{-1/2}(\bar{\mbs{R}} - \mbs{R}) \|^2_F &= \text{tr}\{ (\bar{\mbs{R}}  - \mbs{R}) \mbs{R}^{-1}(\bar{\mbs{R}}  - \mbs{R}) \} \\
&=\text{tr}\{  \bar{\mbs{R}} \mbs{R}^{-1} \bar{\mbs{R}}  \} + \text{tr}\{ \mbs{R}
\} +  \text{const.},
\label{eq:SPICEcovariancematchmulti}
\end{split}
\end{equation}
where $\bar{\mbs{R}} \triangleq \frac{1}{T} \sum^T_{t=1}
\mbs{y}_t\mbs{y}^*_t$ and $T$ is the number of snapshots (possibly $T < N$). The derivative of
\eqref{eq:SPICEcovariancematchmulti} w.r.t. $p_k$ is equal to
\begin{equation*}
- \text{tr}\left\{ \bar{\mbs{R}} \mbs{R}^{-1} \frac{\partial
    \mbs{R}}{\partial p_k} \mbs{R}^{-1} \bar{\mbs{R}} \right \} + w_k
= - \| \mbs{a}^*_k \mbs{R}^{-1} \bar{\mbs{R}} \|^2_2 + w_k.
\end{equation*}
Then the \textsc{Spice} algorithms \eqref{eq:SPICEiterationalt1} and \eqref{eq:SPICEiterationalt2} become
\begin{equation}
\hat{p}^{i+1}_k = \hat{p}^i_k  \| \mbs{a}^*_k \what{\mbs{R}}^{-1}_i \bar{\mbs{R}} \|_2 / w^{1/2}_k \quad (\text{\textsc{Spice}}_{\text{a}})
\end{equation}
and, respectively,
\begin{equation}
\hat{p}^{i+1}_k = \hat{p}^i_k  \| \mbs{a}^*_k \hat{\mbs{R}}^{-1}_i \bar{\mbs{R}} \|^2_2 / w_k  \quad (\text{\textsc{Spice}}_{\text{b}}).
\end{equation}
When the number of snapshots $T \geq N$ one may use a modified cost
function, viz. $\| \mbs{R}^{-1/2}(\bar{\mbs{R}} - \mbs{R}) \bar{\mbs{R}}^{-1/2} \|^2_F$, cf. \cite{StoicaEtAl2011_spicearray}.

\section{Lemma proofs}
\label{app:proofs}

\subsection*{Lemma~1}
The inquality in \eqref{eq:lemmaPinv} follows if we can show that
$\what{\mbs{P}}^{-1} \succeq \mbs{A}^*\what{\mbs{R}}^{-1}\mbs{A}$ or,
equivalently, $\mbs{I}_M - \what{\mbs{P}}^{1/2} \mbs{A}^*
(\mbs{A}\what{\mbs{P}}\mbs{A}^*)^{-1} \mbs{A} \what{\mbs{P}}^{1/2}
\succeq \mbs{0}$; however this is obviously true since the left hand side is the orthogonal projection matrix onto the null space of $\mbs{A} \what{\mbs{P}}^{1/2}$.

\subsection*{Lemma~2}
We have that
\begin{equation*}
\begin{split}
c f(\bar{\mbs{p}}) \bigr|_{\bar{\mbs{p}} = c \mbs{p}} &= c\mbs{y}^* (\mbs{A}c\mbs{P} \mbs{A}^*)^{-1} \mbs{y} + c \sum^{M+N}_{k=1} w_k c p_k \\
&= \mbs{y}^* \mbs{R}^{-1} \mbs{y}+ c^2 \sum^{M+N}_{k=1} w_k p_k \\
&= g(\mbs{p})
\end{split}
\end{equation*}
and thus
$$cf(c\hat{\mbs{p}}) = g(\hat{\mbs{p}}) \leq g(\mbs{p}) = c f(c\mbs{p}), \: \forall \{ p_k \geq 0 \}$$
which implies:
$$f(\hat{\bar{\mbs{p}}}) \leq f({\bar{\mbs{p}}}), \forall \{ \bar{p}_k \geq 0 \},$$
and this concludes the proof.

\subsection*{Lemma~3}
A simple calculation yields:
\begin{equation*}
\begin{split}
&(\mbs{y} - \mbs{B}\mbs{x})^*
\mbs{S}^{-1} (\mbs{y} - \mbs{B}\mbs{x})\\
&= \mbs{y}^* \mbs{S}^{-1}
\mbs{y} - \mbs{y}^* \mbs{S}^{-1} \mbs{B} \mbs{x} - \mbs{x}^* \mbs{B}^*
\mbs{S}^{-1} \mbs{y} + \mbs{x}^* \mbs{B}^* \mbs{S}^{-1} \mbs{B}
\mbs{x}
\end{split}
\end{equation*}
Therefore the criterion in \eqref{eq:SPICEalt} can be re-written as:
$$\mbs{x}^* ( \mbs{B}^* \mbs{S}^{-1} \mbs{B} + \mbs{\Pi}^{-1} )
\mbs{x} - \mbs{x}^* \mbs{B}^* \mbs{S}^{-1} \mbs{y} - \mbs{y}^*
\mbs{S}^{-1} \mbs{B} \mbs{x}+ \text{constant}, $$ which yields the minimizer
\begin{equation}
\hat{\mbs{x}} = ( \mbs{B}^* \mbs{S}^{-1} \mbs{B} + \mbs{\Pi}^{-1} )^{-1} \mbs{B}^* \mbs{S}^{-1} \mbs{y}.
\label{eq:SPICEaltsolution2}
\end{equation}

Next, note that
\begin{equation*}
\begin{split}
(\mbs{B}^* \mbs{S}^{-1} \mbs{B} + \mbs{\Pi}^{-1}) \mbs{\Pi} \mbs{B}^*
&= \mbs{B}^* \mbs{S}^{-1} \mbs{B} \mbs{\Pi} \mbs{B}^* + \mbs{B}^*\\
&= \mbs{B}^* \mbs{S}^{-1} ( \mbs{B} \mbs{\Pi}\mbs{B}^* + \mbs{S})
\end{split}
\end{equation*}
or equivalently,
$$(\mbs{B}^* \mbs{S}^{-1} \mbs{B} + \mbs{\Pi}^{-1})^{-1} \mbs{B}^* \mbs{S}^{-1} =  \mbs{\Pi} \mbs{B}^* ( \mbs{B} \mbs{\Pi} \mbs{B}^*  + \mbs{S})^{-1},$$ which means that \eqref{eq:SPICEaltsolution2} can be re-written as in \eqref{eq:SPICEaltsolution1}.

It remains to evaluate the criterion at $\hat{\mbs{x}}$. Because
\begin{equation*}
\begin{split}
\mbs{y}- \mbs{B}\hat{\mbs{x}} &= (\mbs{I}_M - \mbs{B} \mbs{\Pi}
\mbs{B}^* \mbs{R}^{-1})\mbs{y}\\
&= (\mbs{R} - \mbs{B} \mbs{\Pi}
\mbs{B}^*)\mbs{R}^{-1} \mbs{y} \\
&= \mbs{S} \mbs{R}^{-1} \mbs{y}
\end{split}
\end{equation*}
we have
\begin{equation*}
\begin{split}
&(\mbs{y} - \mbs{B}\hat{\mbs{x}})^* \mbs{S}^{-1} (\mbs{y} -
\mbs{B}\hat{\mbs{x}}) + \hat{\mbs{x}}^* \mbs{\Pi}^{-1} \hat{\mbs{x}}
\\
&= \mbs{y}^* \mbs{R}^{-1} \mbs{S} \mbs{R}^{-1} \mbs{y} +
\mbs{y}^*\mbs{R}^{-1} \mbs{B} \mbs{\Pi} \mbs{B}^* \mbs{R}^{-1} \mbs{y}
\\
&= \mbs{y}^* \mbs{R}^{-1} \mbs{y}
\end{split}
\end{equation*}
which concludes the proof.

\section{Spice for identical noise powers}
\label{app:identicalnoise}

In this case the covariance model becomes:
\begin{equation}
\mbs{R} = \mbs{B}\mbs{\Pi}\mbs{B}^* + \sigma^2 \mbs{I}_N
\end{equation}
and \eqref{eq:SPICEalt} becomes
\begin{equation*}
\min_{\mbs{x}} \| \mbs{y} - \mbs{B}\mbs{x} \|^2_2/\sigma^2 +
\sum^{M}_{k=1} |x_k|^2/p_k = \mbs{y}^*\mbs{R}^{-1} \mbs{y}.
\end{equation*}
It follows that the minimizers $\{ p_k \}$ of the \textsc{Spice}
criterion can also be obtained by minimizing the function:
\begin{equation}
\| \mbs{y} - \mbs{B}\mbs{x} \|^2_2/\sigma^2 +
\sum^{M}_{k=1} |x_k|^2/p_k +
\sum^M_{k=1} w_k p_k + \underbrace{\left( \sum^{M+N}_{k=M+1} w_k \right)}_{w^2}\sigma^2.
\label{eq:SPICEidenticalnoisepower}
\end{equation}
Minimization of \eqref{eq:SPICEidenticalnoisepower} w.r.t. $\sigma^2$
and $\{ p_k \}$ gives:
\begin{equation}
\begin{split}
\sigma^2 &= \| \mbs{y} - \mbs{B}\mbs{x} \|_2 / w \\
p_k &= |x_k| / \sqrt{w_k}, \; k=1,\dots,M.
\end{split}
\label{eq:SPICEsolutionidenticalnoisepower}
\end{equation}
Inserting \eqref{eq:SPICEsolutionidenticalnoisepower} in
\eqref{eq:SPICEidenticalnoisepower} yields (to within a multiplicative
factor):
\begin{equation*}
w \| \mbs{y} - \mbs{B}\mbs{x} \|_2 + \| \text{diag}(\sqrt{w_1},
\dots,\sqrt{w_M} ) \mbs{x} \|_1
\end{equation*}
which is the criterion of the square-root \textsc{Lasso} (with weights for the
$\ell_1$-norm  of $\mbs{x}$). The above proof is more direct than the
one in \cite{Babu&StoicaEtAl2013_spicesqrtlasso,RojasEtAl2013_notespice}.

\section{Cyclic minimization interpretation}
\label{app:cyclic}

The gradient approach in Section~\ref{sec:reweightedspice} is simple and quite flexible;
unlike the cyclic minimization approach in
\cite{StoicaEtAl2011_spicespectral,StoicaEtAl2011_spicearray,Stoica&Babu2012_spicelikes},
the gradient approach produced not only the original algorithms but also different
versions of them. However, the gradient approach \emph{cannot} be used
to conclude the monotonic decrease property used in the convergence
analysis in Section~\ref{sec:statisticalinterpretations}. Indeed, while the function
\begin{equation*}
f(\mbs{p}) = \mbs{y}^* \mbs{R}^{-1} \mbs{y} + \sum^{M+N}_{k=1} w_k p_k
\quad  (w_k \text{ given})
\end{equation*}
is convex, the gradient-based algorithms might overshoot the minimum,
and hence they are not guaranteed to monotonically decrease this
function. To prove such a property we need the cyclic minimization
framework.

Let
\begin{equation}
g(\mbs{\beta},\mbs{p}) = \sum^{M+N}_{k=1} \left( \frac{|\beta_k|^2}{p_k} + w_k p_k \right)
\end{equation}
(the augmented function used by this framework). As shown in
\cite{StoicaEtAl2011_spicespectral,StoicaEtAl2011_spicearray,Stoica&Babu2012_spicelikes}
\begin{equation}
\min_{\mbs{\beta}} g(\mbs{\beta},\mbs{p}) =
\mbs{y}^*\mbs{R}^{-1}\mbs{y} + \sum^{M+N}_{k=1} w_k p_k \quad
(\text{s.t. }\mbs{A}\mbs{\beta}=\mbs{y})
\label{eq:gradcyc_constraint}
\end{equation}
and the minimum is attained at
\begin{equation}
\hat{\beta}_k = p_k \mbs{a}^*_k \mbs{R}^{-1}\mbs{y}.
\label{eq:gradcyc_betak}
\end{equation}
To show this result, let
$$\mbs{\beta} = \begin{bmatrix} \mbs{x} \\ \mbs{y} -\mbs{B} \mbs{x} \end{bmatrix} $$
which satisfies the constraint in \eqref{eq:gradcyc_constraint}; then clearly the result is
equivalent to Lemma~3.

It follows from \eqref{eq:gradcyc_constraint} that to get $\mbs{p}$ that minimizes $f(\mbs{p})$ we
can cyclically minimize $g(\mbs{\beta}, \mbs{p})$ w.r.t. $\mbs{\beta}$
and $\mbs{p}$. For given $\mbs{p}$, the minimizing $\mbs{\beta}$ is
given by \eqref{eq:gradcyc_betak}. For a given $\mbs{\beta}$, the minimization of
$g(\mbs{\beta},\mbs{p})$ w.r.t. $\mbs{p}$ yields
\begin{equation}
\hat{\mbs{p}}_k = |\beta_k| / w^{1/2}_k.
\label{eq:gradcyc_pupdate}
\end{equation}
Iteratively, this means (combining \eqref{eq:gradcyc_betak} and \eqref{eq:gradcyc_pupdate} into one equation):
\begin{equation}
\hat{p}^{i+1}_k = \hat{p}^i_k |\mbs{a}^*_k \mbs{R}^{-1}_i \mbs{y}|/w^{1/2}_k
\end{equation}
which is \eqref{eq:SPICEiterationalt1}. Therefore, for
\eqref{eq:SPICEiterationalt1} the monotonic decreasing property
of $f(\mbs{p})$ is guaranteed:
\begin{equation}
f(\hat{\mbs{p}}^i) = g(\hat{\mbs{\beta}}^i, \hat{\mbs{p}}^i) \geq
g(\hat{\mbs{\beta}}^i, \hat{\mbs{p}}^{i+1}) \geq
g(\hat{\mbs{\beta}}^{i+1}, \hat{\mbs{p}}^{i+1}) \geq
f(\hat{\mbs{p}}^{i+1}).
\label{eq:gradcyc_inequalities}
\end{equation}
But does this property hold for \eqref{eq:SPICEiterationalt2} as well? For \eqref{eq:SPICEiterationalt2}, i.e.,
$$\hat{p}^{i+1}_k = \hat{p}^{i}_k |\mbs{a}^*_k \mbs{R}^{-1}_i
\mbs{y}|^2/w_k $$
we have that:
\begin{equation*}
\begin{split}
\frac{|\hat{\beta}^i_k|^2}{\hat{p}^{i+1}_k} + w_k \hat{p}^{i+1}_k &=
\frac{(\hat{p}^i_k)^2 |\mbs{a}^*_k \what{\mbs{R}}^{-1}_i
  \mbs{y}|^2}{\hat{p}^i_k |\mbs{a}^*_k \what{\mbs{R}}^{-1}_i \mbs{y}|^2}w_k +
w_k \frac{\hat{p}^i_k |\mbs{a}^*_k \what{\mbs{R}}^{-1}_i
  \mbs{y}|^2}{w_k} \\
&= w_k \hat{p}^i_k +
\frac{|\hat{\beta}^i_k|}{\hat{p}^i_k}.
\end{split}
\end{equation*}
Hence
$$g(\hat{\mbs{\beta}}^i, \hat{\mbs{p}}^i) = g(\hat{\mbs{\beta}}^i, \hat{\mbs{p}}^{i+1})$$
and the monotonic decrease property holds for \eqref{eq:SPICEiterationalt2} too (owing to the
second inequality in \eqref{eq:gradcyc_inequalities}).

\bibliography{refs_spicenotes}
\bibliographystyle{ieeetr}

\end{document}